\documentclass[11pt]{amsart}
\usepackage[mathscr]{eucal}
\usepackage{amsfonts}
\usepackage{amsmath}
\usepackage{amsthm}
\usepackage{amssymb}
\usepackage{amscd}
\usepackage{latexsym}

\theoremstyle{plain}
\newtheorem{definition}[equation]{Definition}

\newtheorem{corollary}[equation]{Corollary}
\newtheorem{lemma}[equation]{Lemma}
\newtheorem{proposition}[equation]{Proposition}
\newtheorem{theorem}[equation]{Theorem}

\theoremstyle{definition}

\numberwithin{equation}{subsection}
\renewcommand{\mathfrak}{\mathcal}

\begin{document}
\title{Hochschild and cyclic homology of preprojective algebras of ADE quivers}
\author{Pavel Etingof and Ching-Hwa Eu}
\maketitle
\pagestyle{myheadings}
\markboth{Pavel Etingof and Ching-Hwa Eu}{Hochschild and cyclic homology of preprojective algebras of ADE quivers}

\section{Introduction}

In this paper we compute the additive structure of the Hochschild
\linebreak (co)homology and cyclic homology of preprojective algebras 
of ADE quivers over a field of characteristic zero. 
That is, we compute the (co)homology spaces together with the
grading induced by the natural grading on the preprojective
algebra (in which all edges have degree 1). 
We also use the result (for second cohomology) to find the universal
deformation of the preprojective algebra. 

This generalizes the
results of the papers \cite{ES1}, \cite{ES2}, where 
the dimensions of the Hochschild cohomology groups were found 
for type A and partially for type D. 

Our computation is based on the same method that was used by 
the second author in the paper \cite{Eu}, where the same problem
was solved for centrally extended preprojective algebras,
introduced by E. Rains and the first author. Namely, we use the 
periodic (with period 6) Schofield resolution of the algebra, 
and consider the corresponding complex computing the Hochschild
homology. Using this complex, we find the possible range
of degrees in which each particular Hochschild homology and
space can sit. Then we use this information, as well
as the Connes complex for cyclic homology and the formula for the
Euler characteristic of cyclic homology to find the exact dimensions 
of the homogeneous components  of the homology groups. 
Then we show that the same computation actually yields the
Hochschild cohomology spaces as well.

We note that for connected non-Dynkin quivers, the Hochschild (co)homology
and the cyclic homology of the preprojective algebra were
calculated in \cite{CBEG,EG}; in this case, unlike the ADE case, the homological
dimension of the preprojective algebra is 2, so the situation is simpler. 

{\bf Acknowledgments.} P.E. is grateful to V. Ostrik and
for a useful discussion and to K. Erdmann for references.
The work of the authors was  partially supported by the NSF grant
DMS-0504847. 

\section{Preliminaries}
\subsection{Quivers and path algebras}
Let $Q$ be a quiver of ADE type with vertex set $I$ and $|I|=r$. 
We write $a\in Q$ to say that $a$ is an arrow in $Q$. 

We define $Q^*$ to be the quiver obtained from $Q$ by reversing
all of its arrows. We call $\bar Q=Q\cup Q^*$ the \emph{double}
of $Q$. Let $C$ be the adjacency matrix corresponding to the
quiver $\bar Q$. 

The concatenation of arrows generate the \emph{nontrivial
paths} inside the quiver $\bar Q$. We define $e_i$, $i\in I$ to
be the \emph{trivial path} which starts and ends at $i$. The
\emph{path algebra} $P_{\bar Q}=\mathbb{C}\bar Q$ of $\bar Q$ over $\mathbb{C}$ is the $\mathbb{C}$-algebra with basis the paths in $\bar Q$ and the product $xy$ of two paths $x$ and $y$ to be their concatenation if they are compatible and $0$ if not. We define the \emph{Lie bracket} $[x,y]=xy-yx$.

Let $R=\oplus\mathbb{C}e_i$. Then $R$ is a commutative semisimple
algebra, and $P_Q$ is naturally an
$R$-bimodule.

\subsection{Frobenius algebras}
Let $A$ be a finite dimensional unital $\mathbb{C}-$algebra. We
call it Frobenius if there is a linear function $f:A\rightarrow\mathbb{C}$, such that the form $(x,y):=f(xy)$ is nondegenerate, or, equivalently, if there exists an isomorphism $\phi:A\stackrel{\simeq}{\rightarrow}A^*$ of left $A-$modules: given $f$, we can define $\phi(a)(b)=f(ba)$, and given $\phi$, we define $f=\phi(1)$.

If $\tilde f$ is another linear function satisfying the same
properties as $f$ from above, then $\tilde f(x)=f(xa)$ for some
invertible $a\in A$. Indeed, we define the form $\{a,b\}=\tilde f(ab)$. Then $\{-,1\}\in A^*$, so there is an $a\in A$, such that $\phi(a)=\{-,1\}$. Then $\tilde f(x)=\{x,1\}=\phi(a)(x)=f(xa)$.

\subsection{The Nakayama automorphism}
Given a Frobenius algebra $A$ (with a function $f$ inducing a
bilinear form $(-,-)$ from above), the automorphism
$\eta:A\rightarrow A$ defined by the equation $(x,y)=(y,\eta(x))$
is called the \emph{Nakayama automorphism} (corresponding to
$f$). 

We note that the freedom in choosing $f$ implies 
that $\eta$ is uniquely determined up to an inner automorphism. Indeed,
let $\tilde f(x)=f(xa)$ and define the bilinear form $\{a,b\}=\tilde f(ab)$. Then
\begin{align*}
\{x,y\}&=\tilde f(xy)=f(xya)=(x,ya)=(ya,\eta(x))=f(ya\eta(x)a^{-1}a)\\
&=(y,a\eta(x)a^{-1}).
\end{align*}

\subsection{The preprojective algebra}
Given an ADE-quiver $Q$, we define the \emph{preprojective
algebra} $\Pi_Q$ to be the quotient of the path algebra $P_{\bar Q}$ by
the relation $\sum\limits_{a\in Q}[a,a^*]=0$. It is known that 
$\Pi_Q$ is a Frobenius algebra (see e.g. \cite{ES2},\cite{MOV}).
 From now on, we write $A=\Pi_Q$.

\subsection{Graded spaces and Hilbert series}

Let $W=\oplus_{d\geq0}W(d)$ be a $\mathbb Z_+$-graded
vector space, with finite dimensional homogeneous subspaces. 
We denote by $M[n]$ the same space with
grading shifted by $n$. The graded dual space $M^*$ is defined by the
formula $M^*(n)=M(-n)^*$.  

\begin{definition} \textnormal{(The Hilbert series of vector spaces)}\\
We define the \emph{Hilbert series} $h_W(t)$ to be the series
\begin{displaymath}
h_W(t)=\sum\limits_{d=0}^{\infty}\dim W(d)t^d.
\end{displaymath}
\end{definition}

\begin{definition} \textnormal{(The Hilbert series of bimodules)}\\
Let $W=\oplus_{d\geq0}W(d)$ be a $\mathbb{Z_+}$-graded bimodule
over the ring $R$, 
so we can write $W=\oplus W_{i,j}$. We define the 
\emph{Hilbert series} $H_W(t)$ to be a matrix valued series with the entries 
\begin{displaymath}
H_W(t)_{i,j}=\sum\limits_{d=0}^{\infty}\dim\ W(d)_{i,j}t^d.
\end{displaymath}
\end{definition}

\subsection{Root system parameters}

Let $w_0$ be the longest element of the Weyl group $W$ of $Q$. 
Then we define $\nu$ to be the involution of $I$, such
that $w_0(\alpha_i)=-\alpha_{\nu(i)}$ (where $\alpha_i$ is the
simple root corresponding to $i\in I$). It turns out that
$\eta(e_i)= e_{\nu(i)}$ (\cite{S}; see \cite{ES2}).

Let $m_i$, $i=1,...,r$, be the exponents of the root system
attached to $Q$, enumerated in the increasing order. 
Let $h=m_r+1$ be the Coxeter number of $Q$. 

Let $P$ be the permutation matrix corresponding to the involution
$\nu$. Let $r_+=\dim\ker(P-1)$ and $r_-=\dim\ker(P+1)$.
Thus, $r_-$ is half the number of vertices which are not fixed by
$\nu$, and $r_+=r-r_-$. 

\section{The main results}

Let $U$ be a positively graded vector space with Hilbert series 
$h_U(t)=\sum\limits_{i,\,m_i<\frac{h}{2}}t^{2m_i}$. 
Let $Y$ be a vector space with $\dim Y=
r_+-r_--\#\{i:m_i=\frac{h}{2}\}$, and let
$K=\ker(P+1),\,L=\langle e_i|\nu(i)=i\rangle$, so that $\dim K=r_-$, $\dim L=r_+-r_-$
(we agree that the spaces $K,L,Y$ sit in degree zero). 

The main results of this paper are the following theorems.

\begin{theorem}\label{t1}
The Hochschild cohomology spaces of $A$, as graded spaces, are
as follows: 
\begin{align*}
{HH^0}(A)&=U[-2]\oplus L[h-2],\\
{HH^1}(A)&=U[-2],\\
{HH^2}(A)&=K[-2],\\
{HH^3}(A)&=K[-2],\\
{HH^4}(A)&=U^*[-2],\\
{HH^5}(A)&=U^*[-2]\oplus Y^*[-h-2],\\
{HH^6}(A)&=U[-2h-2]\oplus Y[-h-2],
\end{align*}
and ${HH^{6n+i}}(A)={HH^i}(A)[-2nh]\,\forall
i\geq1$.
\end{theorem}

\begin{corollary}\label{center}
The center $Z=HH^0(A)$ of $A$ has Hilbert series
$$
h_Z(t)=\sum\limits_{i,\,m_i<\frac{h}{2}}t^{2m_i-2}+(r_+-r_-)t^{h-2}. 
$$
\end{corollary}

\begin{theorem}\label{t2}
The Hochschild homology spaces of $A$, as graded spaces, are
as follows:
\begin{align*}
{HH_0}(A)=R,\\
{HH_1}(A)=U,\\
{HH_2}(A)=U\oplus Y[h],\\
{HH_3}(A)=U^*[2h]\oplus Y^*[h],\\
{HH_4}(A)=U^*[2h].\\
{HH_5}(A)=K[2h],\\
{HH_6}(A)=K[2h],\\
\end{align*}
and ${HH_{6n+i}}(A)={HH_i}(A)[2nh]\,\forall i\geq1$. 
\end{theorem}

(Note that the equality $HH_0(A)=R$ was established in
\cite{MOV}). 

\begin{theorem}\label{t3}
The reduced cyclic homology spaces of $A$, as graded spaces, are
as follows:
\begin{align*}
{\overline{HC}_0}(A)=0,\\
{\overline{HC}_1}(A)=U,\\
{\overline{HC}_2}(A)=Y^*[h],\\
{\overline{HC}_3}(A)=U^*[2h],\\
{\overline{HC}_4}(A)=0,\\
{\overline{HC}_5}(A)=K[2h],\\
\end{align*}
and ${\overline{HC}_{6n+i}}(A)={\overline{HC}_i}(A)[2nh]\,\forall i\geq 0$. 
\end{theorem}

The rest of the paper is devoted to the proof of Theorems \ref{t1},\ref{t2},\ref{t3}

\section{Hochschild (co)homology and cyclic homology of A}
\subsection{The Schofield resolution of A}
We want to compute the Hochschild (co)homology of $A$, by using
the Schofield resolution, 
described in \cite{S}. 

Define the 
$A-$bimodule $\mathfrak{N}$ obtained from $A$ by twisting the
right action by $\eta$, i.e., 
$\mathfrak{N}=A$ as a vector space, 
and $\forall a,b\in A,x\in\mathfrak{N}:a\cdot x\cdot b=ax\eta(b).$
Introduce the notation $\epsilon_a=1$ if $a\in Q$,
$\epsilon_a=-1$ if $a\in Q^*$. Let $x_i$ be a homogeneous basis
of $A$ and $x_i^*$
the dual basis under the form attached to the Frobenius algebra
$A$. Let $V$ be the bimodule spanned by the edges of $\bar Q$. 
We start with the following exact sequence:
\[
0\rightarrow\mathfrak{N}[h]\stackrel{i}{\rightarrow}A\otimes_R A[2]\stackrel{d_2}{\rightarrow}A\otimes_RV\otimes_RA\stackrel{d_1}{\rightarrow}A\otimes_RA\stackrel{d_0}{\rightarrow}A\rightarrow 0,
\]
where
\begin{align*}
d_0(x\otimes y)&=xy,\\
d_1(x\otimes v\otimes y)&=xv\otimes y-x\otimes vy,\\
d_2(z\otimes t)&=\sum\limits_{a\in\bar Q}\epsilon_aza\otimes a^*\otimes t+\sum\limits_{a\in\bar Q}\epsilon_az\otimes a\otimes a^*t,\\
i(a)&=a\sum x_i\otimes x_i^*.
\end{align*}

Since $\eta^2=1,$ we can make a canonical identification
$A=\mathfrak{N}\otimes_A\mathfrak{N}$ (via $x\mapsto x\otimes
1$), so by tensoring the above exact sequence with
$\mathfrak{N}$, we obtain the exact sequence

\[
0\rightarrow A[2h]\stackrel{d_6}{\rightarrow}A\otimes_R \mathfrak{N}[h+2]\stackrel{d_5}{\rightarrow}A\otimes_RV\otimes_R\mathfrak{N}[h]\stackrel{d_4}{\rightarrow}A\otimes_R\mathfrak{N}[h]\stackrel{j}{\rightarrow}\mathfrak{N}[h]\rightarrow 0,
\]
and by connecting both sequences with $d_3=ij$ and repeating this process, we obtain the Schofield resolution which is periodic with period $6$:
\begin{align*}
\ldots&\rightarrow A\otimes A[2h]\stackrel{d_6}{\rightarrow}A\otimes_R \mathfrak{N}[h+2]\stackrel{d_5}{\rightarrow}A\otimes_RV\otimes_R\mathfrak{N}[h]\stackrel{d_4}{\rightarrow}A\otimes_R\mathfrak{N}[h]\\
&\stackrel{d_3}{\rightarrow}A\otimes_RA[2]\stackrel{d_2}{\rightarrow}A\otimes_RV\otimes_RA\stackrel{d_1}{\rightarrow}A\otimes_RA\stackrel{d_0}{\rightarrow}A\rightarrow0.
\end{align*}
This implies that the Hochschild homology and cohomology of $A$
is periodic with period $6$, in the sense that the shift of the
(co)homological degree by $6$ 
results in the shift of degree by $2h$ (respectively
$-2h$).

\subsection{The Hochschild homology complex}
Let $A^{op}$ be the algebra $A$ with opposite multiplication. We define $A^e=A\otimes_R A^{op}$. Then any $A-$bimodule naturally becomes a left $A^e-$ module (and vice versa).

We make the following identifications (for all integers $m\geq0$):\\
$(A\otimes_RA)\otimes_{A^e}A[2mh]=A^R[2mh]:\,(a\otimes b)\otimes c=bca$,\\
$(A\otimes_RV\otimes_RA)\otimes_{A^e}A[2mh]=(V\otimes_R A)^R[2mh]:\,(a\otimes x\otimes b)\otimes c=-x\otimes bca$,\\
$(A\otimes_RA)\otimes_{A^e}A[2mh+2]=A^R[2mh+2]:\,(a\otimes b)\otimes c=-bca$,\\
$(A\otimes_R\mathfrak{N})\otimes_{A^e}A[(2m+1)h]=\mathfrak{N}^R[(2m+1)h]:\,(a\otimes b)\otimes c=-b\eta(ca)$,\\
$(A\otimes_RV\otimes_R\mathfrak{N})\otimes_{A^e}A[(2m+1)h]=(V\otimes_R A)^R[(2m+1)h]:$\\$(a\otimes x\otimes b)\otimes c=x\otimes b\eta(ca)$,\\
$(A\otimes_R\mathfrak{N})\otimes_{A^e}A[(2m+1)h+2]=\mathfrak{N}^R[(2m+1)h+2]:\,(a\otimes b)\otimes c=b\eta(ca)$.\\

Now, we apply to the Schofield resolution the functor $-\otimes_{A_e}A$ to calculate the Hochschild homology:

\begin{align*}
\ldots&\rightarrow
\underbrace{A^R[2h]}_{=C_6}\stackrel{d_6'}{\rightarrow}\underbrace{\mathfrak{N}^R[h+2]}_{=C_5}\stackrel{d_5'}{\rightarrow}\underbrace{(V\otimes_R\mathfrak{N})^R[h]}_{=C_4}\stackrel{d_4'}{\rightarrow}\\
&\stackrel{d_4'}{\rightarrow}\underbrace{\mathfrak{N}^R[h]}_{=C_3}\stackrel{d_3'}{\rightarrow}\underbrace{A^R[2]}_{=C_2}\stackrel{d_2'}{\rightarrow}\underbrace{(V\otimes_RA)^R}_{=C_1}\stackrel{d_1'}{\rightarrow}\underbrace{A^R}_{=C_0}\rightarrow0.
\end{align*}

We compute the differentials:
\[
d_1'(a\otimes b)=d_1(-1\otimes a\otimes 1)\otimes_{A^e} b=(-a\otimes 1+1\otimes a)\otimes_{A^e} b=[a,b],
\]
\begin{align*}
d_2'(x)&=d_2(-1\otimes1)\otimes_{A^e} x=-(\sum\limits_{a\in\bar Q}\epsilon_aa\otimes a^*\otimes 1+\sum\limits_{a\in\bar Q}\epsilon_a1\otimes a\otimes a^*)\otimes_{A^e} x\\
&=-\sum\limits_{a\in\bar Q}\epsilon_aa^*\otimes [a,x],
\end{align*}
\begin{eqnarray*}
d_3'(x)&=&d_3(-1\otimes1)\otimes_{A^e}\eta(x)=-\sum (x_i\otimes x_i^*)\otimes_{A^e}\eta(x)=\sum x_i^*\eta(x)x_i\\
&=&\sum x_i^*xx_i=\sum x_ix\eta(x_i^*),
\end{eqnarray*}
the second to last equality is true, since we can assume that each $x_i$ lies in a subspace $e_kAe_{\nu_k}$, and then we see that\\ $x_i^*\eta(x)x_i=x_i^*x_i=x_i^*xx_i$ if $x=e_k$, $k=\nu(k)$,\\and $x_i^*\eta(x)x_i=0=x_i^*xx_i$ if $k\neq\nu(k)$ or $x=e_j$, $j\neq k$ or $\deg x>0$,\\
and the last equality is true because if $(x_i^*)$ is a dual basis of $(x_i)$, then $(x_i)$ is a dual basis of $\eta(x_i^*)$.

\[
d_4'(a\otimes b)=d_4(1\otimes a\otimes 1)\otimes_{A^e} \eta(b)=(a\otimes 1-1\otimes a)\otimes_{A^e} \eta(b)=ab-b\eta(a),
\]
\begin{align*}
d_5'(x)&=d_5(1\otimes1)\otimes_{A^e}x=(\sum\limits_{a\in\bar Q}\epsilon_aa\otimes a^*\otimes 1+\sum\limits_{a\in\bar Q}\epsilon_a1\otimes a\otimes a^*)\otimes_{A^e} x\\
&=\sum\limits_{a\in\bar Q}\epsilon_aa^*\otimes (x\eta(a)-ax),
\end{align*}

\begin{align*}
d_6'(x)&=d_6(1\otimes1)\otimes_{A^e}x=\sum (x_i\otimes x_i^*)\otimes_{A^e}x=
\sum x_i^*\eta(x)\eta(x_i)\\
&=\sum x_i\eta(x)x_i^*=\sum x_ixx_i^*,
\end{align*}
the second to last equality is true because if $(x_i^*)$ is a dual basis of $(x_i)$, then $(x_i)$ is a dual basis of $\eta(x_i^*)$, and \\
the last equality is true because for each $j\in I$, $\sum x_ie_jx_i^*=\sum\dim(e_kAe_j)\omega_j$, where we call $\omega_j$ the dual of $e_j$, and $\dim(e_kAe_j)=\dim(e_kAe_{\nu(j)})$ (given a basis in $e_kAe_j$, the involution which reverses all arrows gives us a basis in $e_jAe_k$, its dual basis lies in $e_kAe_{\nu(j)}$).

Since $A=[A,A]+R$ (see \cite{MOV}), $HH_0(A)=R$, and $HH_6(A)$ sits in degree $2h$.

Let us define $\overline{HH_i}(A)=HH_i(A)$ for $i>0$ and $\overline{HH_i}(A)=HH_i(A)/R$ for $i=0$. Then $\overline{HH_0}(A)=0$.

The top degree of $A$ is $h-2$ (since
$h_A(t)=\frac{1+Pt^h}{1-Ct+t^2}$ 
by \cite[2.3.]{MOV}, and $A$ is finite dimensional).
Thus we see immediately from the homology complex that ${HH_1}(A)$ lives in degrees between $1$ and $h-1$, ${HH_2}(A)$ between $2$ and $h$, ${HH_3}(A)$ between $h$ and $2h-2$, ${HH_4}(A)$ between $h+1$ and $2h-1$, ${HH_5}(A)$ between $h+2$ and $2h$ and ${HH_6}(A)$ in degree $2h$.

\subsection{Self-duality of the homology complex}
The nondegenerate form allows us to make identifications 
$A=\mathfrak{N}^*[h-2]$ and $\mathfrak{N}=A^*[h-2]$ via $x\mapsto(-,x)$. 

We can define a nondegenerate form on $V\otimes A$ and $V\otimes\mathfrak{N}$ by
\begin{equation}
(a\otimes x_a,b\otimes x_b)=\delta_{a,b^*}\epsilon_a(x_a,x_b)
\end{equation}
where $a,b\in Q$, and $\delta_{x,y}$ is $1$ if $x=y$ and $0$
else. 
This allows us to make identifications 
$V\otimes_RA=(V\otimes_R\mathfrak{N})^*[h]$ and $V\otimes_R\mathfrak{N}=(V\otimes_RA)^*[h]$.

Let us take the first period of the Hochschild homology complex, i.e. the part involving the first $6$ bimodules:

\[
\underbrace{\mathfrak{N}^R[h+2]}_{=C_5}\stackrel{d_5'}{\rightarrow}\underbrace{(V\otimes_R\mathfrak{N})^R[h]}_{=C_4}\stackrel{d_4'}{\rightarrow}\underbrace{\mathfrak{N}^R[h]}_{=C_3}\stackrel{d_3'}{\rightarrow}\underbrace{A^R[2]}_{=C_2}\stackrel{d_2'}{\rightarrow}\underbrace{(V\otimes_RA)^R}_{=C_1}\stackrel{d_1'}{\rightarrow}\underbrace{A^R}_{=C_0}\rightarrow0.
\]

By dualizing and using the above identifications, we get the dual complex:

\begin{align*}
\stackrel{(d_3')^*}{\leftarrow}\underbrace{(A^R[2])^*}_{=C_3[-2h]}\stackrel{(d_2')^*}{\leftarrow}\underbrace{((V\otimes_RA)^R)^*}_{=C_4[-2h]}\stackrel{(d_1')^*}{\leftarrow}\underbrace{(A^R)^*}_{=C_5[-2h]}&\leftarrow0.\\
\underbrace{(\mathfrak{N}^R[h+2])^*}_{=C_0[-2h]}\stackrel{(d_5')^*}{\leftarrow}\underbrace{((V\otimes_R\mathfrak{N})^R[h])^*}_{=C_1[-2h]}\stackrel{(d_4')^*}{\leftarrow}\underbrace{(\mathfrak{N}^R[h])^*}_{=C_2[-2h]}&\stackrel{(d_3')^*}{\leftarrow}\\
\end{align*}

We see that $C_i^*=C_{5-i}$. We will now prove that, moreover,
$d_i'=\pm(d_{6-i}')^*$, i.e. the homology complex has a self-duality property. 

\begin{proposition}
One has $d_i'=\pm(d_{6-i}')^*$.
\end{proposition}

\begin{proof}
\underline{$(d_1')^*=d_5'$}:

We have
\begin{align*}
(\sum\limits_{a\in\bar Q}(a\otimes x_a),d'_5(y))&=(\sum\limits_{a\in\bar Q}(a\otimes x_a),\sum\limits_{a\in\bar Q}\epsilon_aa^*\otimes(y\eta(a)-ay))\\
&=\sum\limits_{a\in\bar Q}(x_a,y\eta(a)-ay)=(\sum\limits_{a\in\bar Q}[a,x_a],y)\\
&=(d'_1(\sum\limits_{a\in\bar Q}a\otimes x_a),y)
\end{align*}

\underline{$(d_2')^*=-d_4'$}:

We have
\begin{align*}
(x,d'_4(\sum\limits_{a\in\bar Q}a\otimes x_a))&=(x,\sum\limits_{a\in\bar Q}ax_a-x_a\eta(a))=\sum\limits_{a\in\bar Q}(-[a,x],x_a)\\
&=(\sum\limits_{a\in\bar
Q}\epsilon_aa^*\otimes[a,x],\sum\limits_{a\in \bar Q}a\otimes x_a)
=(-d'_2(x),\sum\limits_{a\in\bar Q}a\otimes x_a)
\end{align*}

\underline{$(d_3')^*=d_3'$}:

We have
\[
(x,d'_3(y))=(x,\sum x_iy\eta(x_i^*))=(\sum x_i^*xx_i,y)=(d'_3(x),y).
\]
\end{proof}

\subsection{Cyclic homology}
Now we want to introduce the cyclic homology which will help us
in computing the Hochschild cohomology of $A$. We have the Connes exact sequence
\[
0\rightarrow \overline{HH_0}(A)\stackrel{B_0}{\rightarrow}\overline{HH_1}(A)\stackrel{B_1}{\rightarrow}\overline{HH_2}(A)\stackrel{B_2}{\rightarrow}\overline{HH_3}(A)\stackrel{B_3}{\rightarrow}\overline{HH_4}(A)\rightarrow\ldots
\]
where the $B_i$ are the Connes differentials (see \cite[2.1.7.]{Lo}) and the $B_i$ are all degree-preserving.
We define the \emph{reduced cyclic homology} (see \cite[2.2.13.]{Lo})

\begin{align*}
\overline{HC_i}(A)&=\ker(B_{i+1}:\overline{HH_{i+1}}(A)\rightarrow\overline{HH_{i+2}}(A))\\
&=\mbox{Im}(B_i:\overline{HH_i}(A)\rightarrow\overline{HH_{i+1}}(A)).
\end{align*}

The usual cyclic homology $HC_i(A)$ is related to the reduced one
by the equality $\overline{HC}_i(A)=HC_i(A)$ for $i$ odd 
and $\overline{HC}_0(A)=HC_0(A)/R$ for $i$ even. 

Let $U={HH_1}(A)$. Then by the degree argument and
the injectivity of $B_1$
(which follows from the fact that $\overline{HH}_0(A)=0$), we
have ${HH_2}(A)=U\oplus Y[h]$ where
$Y={HH_2}(A)(h)$ (the degree-$h$-component). 
Using the duality of the Hochschild homology complex, we find 
${HH_4}(A)=U^*[2h]$ and
${HH_3}(A)=U^*[2h]\oplus Y^*[h]$. 
Let us set $K={HH_5}(A)[-2h]$.

So we can rewrite the Connes exact sequence as follows:
$$
\CD
\text{degree}\\
@.                                0                    \\
@.                              @VVV\\
 1\leq\deg\leq h-1   @.  {HH_1}(A)@=   U     @.@.@.  \overline{HC}_0(A)=0\\
@.                              @VB_1VV    @V\sim VV\\
 1\leq\deg\leq h     @.  {HH_2}(A)@=   U     @. \oplus @.    Y[h] @.  \overline{HC}_1(A)=U\\
@.                              @VB_2VV               @.    @.    @V\sim VV\\
 h+1\leq\deg\leq 2h-1  @.  {HH_3}(A)@= U^*[2h] @. \oplus @. Y^*[h] @. \, \overline{HC}_2(A)=Y^*[h] \\
@.                              @VB_3VV     @V\sim VV\\
h+1\leq\deg\leq2h-1  @.  {HH_4}(A)@= U^*[2h] @.@.@.  \overline{HC}_3(A)=U^*[2h]\\
@.                              @VB_4VV       @V0VV \\
2h    @.  {HH_5}(A)@=   K[2h]     @.@.@.  \overline{HC}_4(A)=0\\
@.                              @VB_5VV     @V\sim VV\\
2h                   @.  {HH_6}(A)@=   K[2h]     @.@.@.  \overline{HC}_5(A)=K[2h]\\
@.                              @VB_6VV       @V0VV    \\
2h+1\leq\deg\leq3h-1\, @.  {HH_7}(A)@=U[2h]\\
@.                              @VB_7VV\\
@.                              \vdots
\endCD
$$

From the exactness of the sequence it is clear that $B_2$ and
$B_3$ restrict to an isomorphism on $Y[h]$ 
and $U^*[2h]$ respectively and that $B_4=0$. $B_6=0$ because 
it preserves degrees, so $B_5$ is an isomorphism.

An analogous argument applies to the portion of the Connes
sequence from homological degree $6n+1$ to $6n+6$ for $n>0$. 

Thus we see that the reduced cyclic homology groups $\overline{HC}_i(A)$ live in different degrees:
$\overline{HC}_{6n+1}(A)$ between $2hn+1$ and $2hn+h-1$,
$\overline{HC}_{6n+2}(A)$ in degree $2hn+h$,
$\overline{HC}_{6n+3}(A)$ between $2hn+h+1$ and $2hn+2h-1$, and
$\overline{HC}_{6n+5}(A)$ in degree $2hn+2h$.
So to prove the main results, it is sufficient to determine 
the Hilbert series of the cyclic homology spaces. 

This is done with the help of the following lemma. 

\begin{lemma}
The Euler characteristic of the reduced cyclic homology $\chi_{\overline{HC}(A)}(t)=\sum\ (-1)^ih_{\overline{HC}_i(A)}(t)$ is 
\[
\sum\limits_{k=0}^\infty a_kt^k=\frac{1}{1-t^{2h}}(-\sum t^{2m_i}-r_-t^{2h}+(r_+-r_-)t^h).
\]
\end{lemma}
\begin{proof}
To compute the Euler characteristic, we use the theorem from \cite{EG} that 
\[
\prod\limits_{k=1}^\infty(1-t^k)^{-a_k}=\prod\limits_{s=1}^\infty\det H_A(t^s).
\]
From \cite[Theorem 2.3.]{MOV} we know that \[H_A(t)=(1+Pt^h)(1-Ct+t^2)^{-1}.\]

Since $r=r_++r_-$,
\[
\det(1+Pt^h)=(1+t^h)^{r_+}(1-t^h)^{r_-}.
\]

From \cite[Proof of Theorem 4.1.2.]{Eu} we know that
\[
\prod_{s=1}^\infty \det(1-Ct^s+t^{2s})=\prod\limits_{k=1}^\infty(1-t^{2k})^{-\#\{i:m_i\equiv k\mod h\}}.
\]
So
\begin{align*}
\prod\limits_{k=1}^\infty(1-t^k)^{-a_k}&=\prod\limits_{s=1}^\infty\det H_A(t^s)\\
&=\prod\limits_{s=1}^\infty(1+t^{hs})^{r_+}(1-t^{hs})^{r_-}\det(1-Ct^s+t^{2s})^{-1}\\
&=\frac{\prod\limits_{s\,even}(1-t^{hs})^{r_-}}{\prod\limits_{s\,odd}(1-t^{hs})^{r_+-r_-}}\prod\limits_{k=1}^\infty(1-t^{2k})^{\#\{i:m_i\equiv k\mod h\}}.
\end{align*}

It follows that
\[
\chi_{\overline{HC}(A)}(t)=\sum\limits_{k=0}^\infty
a_kt^k=(1+t^h+t^{2h}+\ldots)(-\sum t^{2m_i}-r_-t^{2h}+(r_+-r_-)t^h).
\]
This implies the lemma. 
\end{proof}

Since all $\overline{HC}_i(A)$ live in different degrees, we can immediately derive their Hilbert series from the Euler characteristic:
\begin{align*}
h_{\overline{HC}_1(A)}(t)&=\sum\limits_{i,\,m_i<\frac{h}{2}}t^{2m_i},\\
h_{\overline{HC}_2(A)}(t)&=(r_+-r_--\#\{i:m_i=\frac{h}{2}\})t^h,\\
h_{\overline{HC}_3(A)}(t)&=\sum\limits_{i,\,m_i>\frac{h}{2}}t^{2m_i},\\
h_{\overline{HC}_5(A)}(t)&=r_-t^{2h}.\\
\end{align*}
It follows that
$h_U(t)=\sum\limits_{i,\,m_i<\frac{h}{2}}t^{2m_i}$, $\dim
Y=r_+-r_--\#\{i:m_i=\frac{h}{2}\}$, $\dim K=r_-$, and $Y,K$ sit
in degree zero. 

This completes the proof of Theorems \ref{t2},\ref{t3}.
 
\subsection{The Hochschild cohomology complex}
Now we would like to prove Theorem \ref{t1}. 

We make the following identifications:
$Hom_{A^e}(A\otimes_R A,A)=A^R$ and $Hom_{A^e}(A\otimes_R \mathfrak{N},A)=\mathfrak{N}^R$,
by identifying $\phi$ with the image $\phi(1\otimes1)=a$ (we write $\phi=a\circ -$), and
$Hom_{A^e}(A\otimes_R V\otimes_R A,A)=(V\otimes_R A)^R[-2]$ 
and $Hom_{A^e}(A\otimes_R V\otimes_R \mathfrak{N},A)=(V\otimes_R \mathfrak{N})^R[-2],$
by identifying $\phi$ which maps $1\otimes a\otimes1\mapsto x_a$
($a\in\bar Q$) with the element $\sum\limits_{a\in\bar Q}\epsilon_{a^*}a^*\otimes x_a$ (we write $\phi=\sum\limits_{a\in\bar Q}\epsilon_{a^*}a^*\otimes x_a\circ -$).

Now, apply the functor $Hom_{A^e}(-,A)$ to the Schofield resolution to obtain the Hochschild cohomology complex

\begin{align*}
\stackrel{d_4^*}{\leftarrow}\mathfrak{N}^R[-h]\stackrel{d_3^*}{\leftarrow}A^R[-2]\stackrel{d_2^*}{\leftarrow}(V\otimes A)^R[-2]\stackrel{d_1^*}{\leftarrow}A^R&\leftarrow0
\\
\ldots\leftarrow A^R[-2h]\stackrel{d_6^*}{\leftarrow}\mathfrak{N}^R[-h-2]\stackrel{d_5^*}{\leftarrow}(V\otimes\mathfrak{N})^R[-h-2]&\stackrel{d_4^*}{\leftarrow}
\end{align*}

\begin{proposition}
Using the differentials $d_i'$ from the Hochschild homology complex, we can rewrite the Hochschild cohomology complex in the following way:
\begin{align*}
\stackrel{d_5'[-2h-2]}{\leftarrow}\mathfrak{N}^R[-h]\stackrel{d_6'[-2h-2]}{\leftarrow}A^R[-2]\stackrel{d_1'[-2]}{\leftarrow}(V\otimes A)^R[-2]\stackrel{d_2'[-2]}{\leftarrow}A^R&\leftarrow0\\
\ldots\leftarrow A^R[-2h]\stackrel{d_3'[-2h-2]}{\leftarrow}\mathfrak{N}^R[-h-2]\stackrel{d_4'[-2h-2]}{\leftarrow}(V\otimes\mathfrak{N})^R[-h-2]&\stackrel{d_5'[-2h-2]}{\leftarrow}.
\end{align*}
\end{proposition}
\begin{proof}
\[
d_1^*(x)(1\otimes a\otimes1)=x\circ d_1(1\otimes a\otimes1)=x\circ(a\otimes1-1\otimes a)=[a,x],
\]
so \[d_1^*(x)=\sum\limits_{a\in\bar Q}\epsilon_{a^*}a^*\otimes[a,x]=d_2'(x).\]

\begin{align*}
d_2^*(\sum\limits_{a\in\bar Q}a\otimes x_a)(1\otimes1)&=(\sum\limits_{a\in\bar Q}a\otimes x_a)\circ(\sum\limits_{b\in\bar Q}\epsilon_bb\otimes b^*\otimes1+\sum\limits_{b\in\bar Q}\epsilon_b1\otimes b\otimes b^*)\\
&=\sum\limits_{a\in\bar Q}(ax_a-x_aa)=\sum\limits_{a\in\bar Q}[a,x_a],
\end{align*}
so
\[
d_2^*(\sum\limits_{a\in\bar Q}a\otimes x_a)=\sum\limits_{a\in\bar Q}[a,x_a]=d_1'(\sum\limits_{a\in\bar Q}a\otimes x_a).
\]

\[
d_3^*(x)(1\otimes1)=x\circ d_3(1\otimes1)=x\circ(\sum x_i\otimes x_i^*)=\sum x_ix x_i^*,
\]
so
\[
d_3^*(x)=\sum x_ixx_i^*=d_6'(x).
\]

\[
d_4^*(x)(1\otimes a\otimes1)=x\circ d_1(1\otimes a\otimes1)=x\circ(a\otimes1-1\otimes a)=ax-x\eta(a),
\]
so \[d_4^*(x)=\sum\limits_{a\in\bar Q}\epsilon_{a^*}a^*\otimes (ax-x\eta(a))=d_5'(x).\]

\begin{align*}
d_5^*(\sum\limits_{a\in\bar Q}a\otimes x_a)(1\otimes1)&=(\sum\limits_{a\in\bar Q}a\otimes x_a)\circ(\sum\limits_{b\in\bar Q}\epsilon_bb\otimes b^*\otimes1+\sum\limits_{b\in\bar Q}\epsilon_b1\otimes b\otimes b^*)\\
&=\sum\limits_{a\in\bar Q}(ax_a-x_a\eta(a)),
\end{align*}
so
\[
d_5^*(\sum\limits_{a\in\bar Q}a\otimes x_a)=\sum\limits_{a\in\bar Q}(ax_a-x_a\eta(a))=d_4'(\sum\limits_{a\in\bar Q}a\otimes x_a).
\]

\[
d_6^*(x)(1\otimes1)=x\circ d_6(1\otimes1)=x\circ(\sum x_i\otimes x_i^*)=\sum x_ix\eta(x_i^*),
\]
so
\[
d_6^*(x)=\sum x_ix\eta(x_i^*)=d_3'(x).
\]
\end{proof}

Thus we see that each 3-term portion of the cohomology complex
can be identified, up to shift in degree, with an appropriate
portion of the homology complex. 

This fact, together with Theorem \ref{t2}, implies Theorem
\ref{t1}.

\section{The deformed preprojective algebra} 

In this subsection we would like to consider the universal
deformation of the preprojective algebra $A$. 
If $\nu=1$, then $P=1$ and hence by Theorem \ref{t1} 
$HH^2(A)=0$ and thus $A$ is rigid. On the other hand, 
if $\nu\ne 1$ (i.e. for types $A_n$, $n\ge 2$, $D_{2n+1}$, and $E_6$),
then $HH^2(A)$ is the space $K$ of $\nu$-antiinvariant functions 
on $I$, sitting in degree $-2$. 

\begin{proposition}\label{uni} Let $\lambda$ be a weight (i.e. a complex
function on $I$) such that $\nu\lambda=-\lambda$. Let $A_\lambda$
be the quotient of $P_Q$ by the relation 
$$
\sum_{a\in Q}[a,a^*]=\sum \lambda_i e_i. 
$$
Then ${\rm gr}A_\lambda=A$ (under the filtration by length of
paths). Moreover, $A_\lambda$, with $\lambda$ a formal parameter 
in $K$, is a universal deformation of $A$. 
\end{proposition}

\begin{proof}
To prove the first statement, it is sufficient to show that for
generic $\lambda$ such that $\nu(\lambda)=-\lambda$, the
dimension of the algebra $A_\lambda$ is the same as the dimension
of $A$, i.e. $rh(h+1)/6$. But by Theorem 7.3 of \cite{CBH}, 
$A_\lambda$ is Morita equivalent to the preprojective algebra of
a subquiver $Q'$ of $Q$, and the dimension vectors 
of simple modules over $A_\lambda$ are known (also from
\cite{CBH}). This allows one to compute the dimension of
$A_\lambda$ for any $\lambda$, and after a somewhat tedious 
case-by-case computation one finds that indeed $\dim
A_\lambda=\dim A$ for a generic $\lambda\in K$. 

The second statement boils down to the fact that the induced 
map $\phi: K\to HH^2(A)$ defined by the above deformation 
is an isomorphism (in fact, the identity). This is 
proved similarly to the case of centrally
extended preprojective algebras, which is considered in \cite{Eu}.
\end{proof}

{\bf Remark.} For type $A_n$ (but not $D$ and $E$) the algebra
$A_\lambda$ for generic $\lambda\in K$ is actually semisimple, 
with simple modules of dimensions $n,n-2,n-4...$.

\end{document}